\documentclass{amsart}
\usepackage{graphicx}
\usepackage{amssymb}
\usepackage{epstopdf}
\usepackage{nicefrac}
\usepackage{enumerate}
\usepackage{hyperref}
\usepackage{float}
\usepackage{color}
\usepackage{pst-all}
\usepackage{tabularx}
\input xy
\xyoption {all}

\newtheorem{thm}{THEOREM}[section]

\newtheorem{defn}[thm]{DEFINITION}

\newtheorem{prob}[thm]{PROBLEM}
\newtheorem{prop}[thm]{PROPOSITION}

\newcommand{\ds}{\displaystyle}

\newcommand{\F}{{\mathcal F}}

\newcommand{\G}{\Gamma}

\newcommand{\e}{{\epsilon}}

\newcommand{\mC}{{\mathbb C}}

\newcommand{\mQ}{{\mathbb Q}}
\newcommand{\mP}{{\mathbb P}}
\newcommand{\mR}{{\mathbb R}}
\newcommand{\mS}{{\mathbb S}}
\newcommand{\mT}{{\mathbb T}}

\newcommand{\mZ}{{\mathbb Z}}

\newcommand{\cP}{{\mathcal P}}

\newcommand{\cR}{{\mathcal R}}

\newcommand{\wtnu}{{\widetilde{\nu}}}

\newcommand{\fV}{{\mathfrak{V}}}

\newcommand{\GF}{\Gamma_{\F}}

\newcommand{\oG}{\overline{\Gamma}}

\newcommand{\ocR}{\overline{\mathcal R}}

\begin{document}

\title{Rigid secondary characteristic classes} 

\author{Steven Hurder}
\address{Steven Hurder, Department of Mathematics, University of Illinois at Chicago, 322 SEO (m/c 249), 851 S. Morgan Street, Chicago, IL 60607-7045}
\email{hurder@uic.edu}
\thanks{Version date: March 14, 2024}

\date{}


  
  \begin{abstract}
We construct families of non-trivial universal rigid secondary classes for   foliations, and then discuss their application to prove that foliations are not homotopic. An observation of Lawson about the non-triviality of the normal Pontrjagin classes of foliations is extended, and then used to construct new families of examples of foliations with non-trivial rigid secondary classes. Examples are given of (abstractly constructed) foliations of compact manifolds with homotopic tangent bundles, but which are not homotopic as foliations.
\end{abstract}

\maketitle

\section{Introduction} \label{sec-intro}

A basic problem in foliation theory is whether a given smooth manifold $M$ has at least one foliation of codimension $q$. The survey of foliation theory  by Lawson \cite{Lawson1974} discusses this problem at length.  An alternate  formulation of  this problem  is to introduce the space 
  $\F_q(M)$   of   foliations of   $M$ with codimension $q$,  with the topology induced from the Grassmannian space of subbundles of $TM$ of codimension $q$, and consider  whether the set $\F_q(M)$ is non-empty, and if so, what are the topological properties of this space.
 
 Lawson   discusses various notions of equivalence between foliations of a   manifold $M$  in \cite[Section~5]{Lawson1974}, and in particular that of homotopic foliations.
Given foliations $\F_0 , \F_1 \in \F_q(M)$ we say they lie in the same smooth path component of $\F_q(M)$  if there exists a smooth path $\F_t \in \F_q(M)$ of foliations between them. That is, the collection $\{\F_t \mid 0 \leq t \leq 1\}$ defines a smooth foliation of codimension $q+1$ on $M \times [0,1]$. We then say that $\F_0$ and $\F_1$ are smoothly homotopic.
\begin{prob}\label{prob-pi0}
Determine the set $\pi_0^{\infty}(\F_q(M))$ of smooth path components of $\F_q(M)$.
\end{prob}
When $M$ has dimension $m$ and $q=(m-1)$, so the leaves of $\F$ are $1$-dimensional, then the set  $\pi_0^{\infty}(\F_q(M))$ equals  the set of  homotopy classes of $1$-dimensional subbundles of $TM$. On the other hand, when $1 \leq q < m-1$, the  space $\F_q(M)$ is far more mysterious, and there is no direct approach to understanding its set of smooth path components. Also, in the definition of homotopy it becomes  important whether the induced foliation on $M \times [0,1]$ is smooth, or just continuous, as discussed by Rosenberg and Thurston \cite{RT1973}. In any case, if $\F_0$ and $\F_1$ are homotopic foliations, then their tangent and normal bundles are homotopic as subbundles of $TM$, and Lawson proposed:
\begin{prob}\cite[Problem~3]{Lawson1974} \label{prob-homotopy}
Show  there   foliations of a smooth manifold $M$ of codimension $q$ which have homotopic tangential distributions, but are not foliated homotopic.
\end{prob}
In this work, we use the non-trivial rigid secondary classes to distinguish classes in  $\pi_0^{\infty}(\F_q(M))$ for a variety of compact manifolds $M$, and give examples of solutions to Problem~\ref{prob-homotopy}.

Recall that the secondary characteristic classes of smooth foliations (see Section~\ref{sec-secondary})   distinguish   foliations up to some notion of equivalence, such as diffeomorphism or foliated concordance. Examples of Thurston \cite{Thurston1972} show that the Godbillon-Vey class for codimension-one foliations can vary non-trivially along a smooth path of foliations. The examples of Heitsch \cite{Heitsch1978} show that the generalized Godbillon-Vey classes for codimension $q > 1$   also vary along smooth paths  of foliations.  The works of Heitsch    in \cite{Heitsch1973,Heitsch1975} and Gel{'}fand,  Fe\u{\i}gin,    and Fuks in  \cite{GFF1974}, showed that this variation can be expressed in terms of a foliated cohomology space associated to $\F$. 

\eject

In contrast,    there are a subset of the secondary classes,   called \emph{rigid classes}, that depend  only on the smooth homotopy   class of the foliation. These classes are best understood in terms of   the classifying spaces of foliations. 
Let $B\G_q$ denote Haefliger's classifying space of smooth foliations of codimension $q$, as introduced in \cite{Haefliger1970,Haefliger1971}; see also  Lawson \cite{Lawson1977}.  
  A smooth foliation $\F$ on $M$ induces a classifying map $\wtnu_{\F} \colon M \to B\G_q$. If $\F_0$ and $\F_1$ are  foliations of $M$ and their classifying maps $\wtnu_{\F_0}$ and $\wtnu_{\F_1}$ are homotopic, then there is a  $\G_q$-structure on  $M \times [0,1]$ which restricts to the foliations $\F_0$ on $M \times \{0\}$ and $\F_1$ on $M \times \{1\}$. 

Now suppose that     $\{\F_t \mid 0 \leq t \leq 1\}$ is a smooth homotopy between $\F_0$ and $\F_1$, then there is an induced classifying map $\wtnu_{\F_t} \colon M \times [0,1] \to B\G_{q+1}$.
Let  $\times_{\mR} \colon B\G_q \to B\G_{q+1}$ be the map defined by taking the product of a codimension-q foliated manifold with $\mR$, where the factor $\mR$ is considered to be transverse to the leaves. 
    The classifying map $\wtnu_{\F_t}$ is a homotopy between 
    the compositions $\times_{\mR} \circ \wtnu_{\F_0} \colon M \to B\G_{q+1}$ and $\times_{\mR} \circ \wtnu_{\F_1} \colon M \to B\G_{q+1}$.  The image of $\times_{\mR}^* \colon H^*(B\G_{q+1} ; \mR) \to H^*(B\G_q ; \mR)$ yields invariants of the smooth homotopy class of a foliation.

 Let $Q$ be the normal bundle to a foliation $\F$ on $M$, then there is a classifying map $\nu_{\F} \colon M \to BO(q)$, such that $Q$ is the pull-back via $\nu_{\F}$ of the universal $\mR^q$-bundle over $BO(q)$. There is also a classifying map $\nu \colon B\G_q \to BO(q)$ for the  normal bundle to the $\G_q$-structure on $B\G_q$.  Let    $B\oG_q$  denote the homotopy fiber of the map $\nu \colon B\G_q \to BO(q)$, then the space  $B\oG_q$ classifies the smooth foliations of codimension $q$ with \emph{framed} normal bundles.  There is a natural map $\overline{\times}_{\mathbb R} \colon B\oG_{q} \to B\oG_{q+1}$, where a framing of the normal bundle to $\F$ is extended by the tangent vector $\partial/\partial t$ to the ${\mathbb R}$-factor. 

Suppose there is a smooth homotopy $\{\F_t \mid 0 \leq t \leq 1\}$   between $\F_0$ and $\F_1$, and $s_0$ is a framing of  the normal bundle to $\F_0$ then  parallel transport along the  second factor in $M \times [0,1]$ induces a framing $s_1$    for the normal bundle to $\F_1$. Let $\wtnu_{\F_i}^{s_i} \colon M \to B\oG_{q}$ denote the classifying map for the   foliation $(\F_i, s_i)$ with framed normal bundle, for $i=0,1$.
Then as before,  there is a homotopy   between 
    the compositions $\overline{\times}_{\mR} \circ \wtnu_{\F_0}^{s_0} \colon M \to B\oG_{q+1}$ and $\overline{\times}_{\mR} \circ \wtnu_{\F_1}^{s_1} \colon M \to B\oG_{q+1}$.  
    It follows that classes in the image of $\overline{\times}_{\mR}^* \colon H^*(B\oG_{q+1} ; \mR) \to H^*(B\oG_q ; \mR)$ are invariants of the smooth homotopy class of a framed foliation $(\F, s)$.
Introduce the subspaces
 \begin{eqnarray}
 \cR_q & = & {\rm Image}\{ \times_{\mR}^* \colon H^*(B\G_{q+1} ; \mR) \to H^*(B\G_q ; \mR)   \} \label{eq-unframedrigid}\\
 \ocR_q & = &  {\rm Image}\{ \times_{\mR}^* \colon H^*(B\oG_{q+1} ; \mR) \to H^*(B\oG_q ; \mR)  \}  \label{eq-framedrigid} \ .
 \end{eqnarray}
The discussion above then yields:
  \begin{thm}\label{thm-nontrivial-rigid2}
  The classes in $\cR_q$ are  invariants of the smooth homotopy class of a foliation; those in 
  $\ocR_q$ are invariants of the smooth homotopy class  of a foliation  with framed normal bundle.
\end{thm}
The only known non-trivial elements of $\cR_q$ are given by the Pontrjagin classes to the normal bundles of foliations.
 The fiber $B\oG_q$ is $(q+1)$-connected  \cite{Thurston1974}, so the map $\nu^* \colon H^{4k}(B{\rm O}(q) ; \mR) \to H^{4k}(B\G_q  ; \mR)$ is injective for $4k \leq (q+2)$, and the Bott Vanishing Theorem \cite{Bott1970,Bott1971} states that $\nu^*$ vanishes for $4k > 2q$.
Define the subspaces of $\cR_q$
\begin{eqnarray}
\cP_{q}^- & = & {\rm Image}\{\nu^* \colon H^*(B{\rm O}(q) ; \mR) \to H^*(B\G_q  ; \mR)  \mid 4\leq  * \leq  q+2\} \\
\cP_{q}^+ & = &  {\rm Image}\{\nu^* \colon H^*(B{\rm O}(q) ; \mR) \to H^*(B\G_q  ; \mR)  \mid q+2 < * \leq 2q\}   \ .   
 \end{eqnarray}
Thus, we have 
 $\cP_q^- \cong \{H^{4k}(B{\rm O}(q) ; \mR) \mid 4 \leq 4k \leq q+2\}$.  The non-triviality of classes in $\cP_q$ is discussed in Section~\ref{sec-pontrjagin}, and only partial results are known. In either case, the  classes in $\cP_{q}^-$ and $\cP_{q}^+$ depend only on the homotopy class of the normal bundle to a foliation, so are not useful for  the study of Problem~\ref{prob-homotopy}.
 On the other hand,   non-trivial classes in $\overline{\cR}_q$ can be used to construct solutions to Problem~\ref{prob-homotopy}. 
 
Our first result is the construction of families of  linearly independent secondary classes in $\ocR_q$ which is done using a   recursive procedure, and described in detail in Section~\ref{sec-secondary}.   
\begin{thm}\label{thm-main2}
For   $q = 2k$ with $k\geq 2$, the subspace $\ocR_q$  is non-trivial, and its  dimension   $\dim(\ocR_q)$ tends to infinity quadratically in $k$, as   $k$ tends to infinity. 
\end{thm}
 The conclusions of Theorem~\ref{thm-main2} are of interest with regards the study of the cohomology of classifying spaces for foliations, but for applications these results are not particularly useful, as one is interested in the set of homotopy classes $[M, B\oG_q]$ of maps of a connected manifold $M$ into $B\oG_q$. For applications, the following notion plays a fundamental role.
 
\begin{defn}\label{def-spherical}
Let $Y$ be a connected   space.  The degree $n$ \emph{spherical cohomology} of $Y$ is:
\begin{equation}
H^n_s(Y ; \mR) = \{\omega \in H^n(Y; \mR) \mid \exists  \ z \in \pi_n(Y) ~ {\rm such \ that} ~ \langle \omega, h(z) \rangle \ne 0\} \ ,
\end{equation}
where $h \colon \pi_n(Y) \to H_n(Y ; \mZ)$ is the Hurewicz homomorphism.
\end{defn}

In Section~\ref{sec-spherical} we give a collection of techniques for showing that a large subset $\ocR^s_q \subset \ocR_q$, as defined in \eqref{eq-4k-rigid} and  \eqref{eq-4k-2-rigid},  of the linearly independent classes in Theorem~\ref{thm-main2} are linearly independent when restricted to the image of the Hurewicz map.  The main conclusion of this section is then:

 \begin{thm}\label{thm-rigidspherical}
For all even $q \geq 4$, the collection $\Delta(\ocR_q^s) \subset H^*_s(B\oG_q ; \mR)$ is a non-empty,  linearly independent  set of spherically supported rigid secondary classes.
 \end{thm}

Finally, in Section~\ref{sec-applications} we give applications  of Theorem~\ref{thm-rigidspherical} to construct infinite families of foliations on compact manifolds which are non-homotopic, but have homotopic tangent bundles, thus are solutions to Problem~\ref{prob-homotopy}.  

 There is an alternate approach to showing  there are non-trivial solutions to Problem~\ref{prob-homotopy}, which uses the dual homotopy invariants for foliations constructed in \cite{Hurder1981a}.   This is the content of    \cite[Theorem~1]{Hurder1981b}  and  \cite[Corollary~4.8]{Hurder1985a}. The purpose of this paper is to explain these results, and state them more precisely, using more elementary arguments without requiring the use of minimal model theory.

\section{Pontrjagin classes} \label{sec-pontrjagin}

Let $B\G_q^+$ denote the classifying space for smooth foliations of codimension $q$ with orientable normal bundles. There is a classifying map $\nu_q \colon B\G_q^+ \to B{\rm SO}(q)$ of the normal bundle to the universal $\G_q^+$-structure on $B\G_q^+$.
In this section, we analyze the induced map $\nu_q^* \colon H^*(B{\rm SO}(q) ; \mQ) \to H^*(B\G_q^+ ; \mQ)$. 
 
 Recall from the definition of the homotopy fiber $B\oG_q$ that there is a sequence of fibrations,
 \begin{equation}\label{eq-lhes}
{\rm SO}(q) \cong \Omega B{\rm SO}(q) \stackrel{\beta}{\longrightarrow} B\oG_q \stackrel{\iota}{\longrightarrow} B\G_q^+ \stackrel{\nu_q}{\longrightarrow} B{\rm SO}(q) \ .
\end{equation}
Some facts are known about  each of these maps $\beta, \iota, \nu_q$, but in general their study presents many open problems.
The seminal observation was made by Bott around 1970. The cohomology ring $H^*(BO(q); \mR) \cong \mR[p_1 , \ldots , p_k]$ where $2k \leq q$, and $p_{j}$ has graded degree $4 j$.

 \begin{thm} \cite{Bott1970}  \label{BVT} Let $\F$ be a codimension-$q$, $C^2$-foliation, and  $\nu_{Q}   \colon M \to BO(q)$   the classifying map for the normal bundle $Q$. Then $\nu_Q^* \colon H^{\ell}(BO(q); \mQ) \to H^{\ell}(M; \mQ)$ is the trivial map for $\ell > 2q$.
\end{thm}

  Thurston showed in \cite{Thurston1974} that the homotopy fiber $B\oG_q$ of $\nu_q$ is $(q+1)$-connected. Thus, the universal map  $\nu_q^* \colon H^{\ell}(BO(q); \mQ) \to H^{\ell}(B\G_q^+; \mQ)$  is injective for degrees up to $(q+2)$. 
  
The group ${\rm SO}(q)$ acts on $B\oG_q$ by twisting the framing of the universal normal bundle to the $\G_q$-structure on $B\oG_q$. This action was studied in the author's works \cite{Hurder1981b,Hurder1985a,Hurder1993}.

First,  recall some elementary facts about ${\rm SO}(q)$ and $B{\rm SO}(q)$ (see \cite[Chapter 15]{MilnorStasheff1974}). 
The key point is there is  a fibration 
  ${\rm SO}(q) \to {\rm SO}(q+1) \stackrel{\pi}{\longrightarrow} \mS^q$ for all $q \geq 1$.
Let $\chi_q \in H^q(\mS^q ; \mR)$ be the integral class which generates the space.

  For $q=1$ we have ${\rm SO}(1)$ is a point, so $\pi \colon {\rm SO}(2) \to \mS^1$ is a homeomorphism. For $q=2$ we obtain the Hopf fibration ${\rm SO}(2) \to {\rm SO}(3) \to \mS^2$ so that ${\rm SO}(3) \cong \mS^3$, and $H^3({\rm SO}(3) ; \mR)$ is generated by the class $T\chi_2 \wedge \pi^*(\chi_2)$.
  For even $q > 2$, we have $H^{2q-1}({\rm SO}(q+1) ; \mR)$ is generated by 
  $T\chi_q \wedge \pi^*\chi_q$.  Moreover,  $T\chi_q \wedge \pi^*\chi_q$ is non-zero on the image of the Hurewicz map $h \colon \mZ \cong \pi_{2q-1}({\rm SO}(q)) \to H_{2q-1}({\rm SO}(q) ; \mZ)$.
   For odd $q \geq 3$, we have $\pi^*(\chi_q) \in H^{q}({\rm SO}(q) ; \mR)$ is non-zero, and   is non-zero on the image of the Hurewicz map $h \colon \mZ \cong \pi_{q}({\rm SO}(q)) \to H_{q}({\rm SO}(q) ; \mZ)$.
\eject

By the homotopy equivalence $\Omega B{\rm SO}(q) \cong {\rm SO}(q)$ we conclude that for $n > 1$:
\begin{itemize}
\item for $q = 2n$, $H^*(B{\rm SO}(q) ; \mQ) \cong \mQ[p_1, \ldots, p_{n-1}, e_{q}]$
\item for $q = 2n+1$, $H^*(B{\rm SO}(q) ; \mQ) \cong \mQ[p_1, \ldots, p_{n}]$.
\end{itemize}
Here, $p_i$ denotes the Pontrjagin class of degree $4i$, and $e_q$ is the Euler class of degree $q$. 
Each class $p_i$ pairs non-trivially with a class in the image of the Hurewicz map $h \colon \pi_{4i}(B{\rm SO}(q)) \to H_{4i}(B{\rm SO}(q) ; \mZ)$.
 
 Let $\xi \to X$ be a real vector bundle of dimension $q$ over a space $X$, and let    $f_{\xi} \colon X \to B{\rm SO}(q)$ be the classifying map for $\xi$, which is unique up to homotopy. Then $p_i(\xi) = f_{\xi}^*(p_i) \in H^{4i}(X; \mQ)$ is the $i^{th}$ rational Pontrjagin class of $\xi$, and the   total   Pontrjagin class of a bundle $\xi$   is given by
\begin{equation}
p(\xi) = 1 + p_1(\xi) + \cdots + p_k(\xi)
\end{equation}
where $k =  \lfloor q/2 \rfloor$ is the largest integer less than or equal to $q/2$.
Given real vector bundles $\xi \to X$ and $\eta \to X$, their total \emph{rational} Pontrjagin classes satisfy the multiplicative formula
\begin{equation}\label{eq-multiplicative}
p(\xi \oplus \eta)  = p(\xi) \ p(\eta) \in H^*(X   ; \mQ) \ .
\end{equation}
In particular, for each $k \geq 1$ we have the formula
\begin{equation}\label{eq-pk}
p_k(\xi \oplus \eta) = \sum_{i+j = k} \ p_i(\xi) \cdot p_j(\eta) \ , 
\end{equation}

 Next, let $X(q)$ denote the $(q+2)$-subcomplex of the standard CW decomposition of $B{\rm SO}(q)$. The fiber $B\oG_q$ of the map $\nu_q$ is $(q+1)$-connected, so there exists a section $\iota_q \colon X(q) \to B\G_q^+$ of $\nu_q$.

 The cohomology ring $H^*(X(q); \mQ)$ is isomorphic to the polynomial ring $H^*(B{\rm SO}(q) ; \mQ)$ truncated in degrees greater than $q+2$.
   For example, $X(2)$ is just the complex projective space $\mC\mP_2$ of dimension $2$ whose cohomology ring is $H^*(X(2) ;  \mR) \cong \mR[e_2]/e_2^3$, and  
    $H^*(X(3) ; \mQ)  \cong \mQ[p_1]/p_1^2$.
      More generally, for $k \geq 1$ we have:
 \begin{eqnarray*}
 H^*(X(4k) ; \mQ) ~ & \cong & ~   \mQ[p_1, \ldots, p_k , e_{4k}] /  \{* > 4k+2\} \\
 H^*(X(4k+1) ; \mQ)  ~ & \cong & ~  \mQ[p_1, \ldots, p_{k}] /  \{* > 4k+3\} \\
 H^*(X(4k+2) ; \mQ)  ~ & \cong & ~  \mQ[p_{1}, \ldots, p_{k+1}, e_{4k+2}] /  \{* > 4k+4\} \\
H^*(X(4k +3) ; \mQ)  ~ & \cong & ~ \mQ[p_1, \ldots , p_{k+1}] / \{* > 4k+5\} \ .
\end{eqnarray*}

 We now apply these observations to show the non-triviality   of   $\nu_q^* \colon H^*(B{\rm SO}(q) ; \mQ) \to H^*(B\G_q^+ ; \mQ)$ 
 for the range of degrees $q+2 < * \leq 2q$.

It was remarked by Lawson    \cite[page  397]{Lawson1974} that the Bott Vanishing Theorem  is optimal, in that the powers of the first Pontrjagin class $p_1^k \in H^{4k}(B\G_{2k} ; \mR)$ are non-vanishing for all $k \geq 1$. This follows from a simple observation. Let $q = q_1 + \cdots + q_{\ell}$ be a partition of $q$, where each $q_i \geq 1$. Then there is a commutative diagram:
 \begin{align}\label{eq-commutingdiagram}
\xymatrix{
& B\G_{q_1}^+ \times \cdots \times B\G_{q_{\ell}}^+  \ar[d]^{\nu_{q_1} \times \cdots \times \nu_{q_{\ell}}}   \ar[r]
& \hspace{2mm} B\Gamma_q^+  \ar[d]^{\nu_q}  \\
X(q_1) \times \cdots \times X(q_{\ell}) \ar[r] \ar[ur]^{\iota_{q_1} \times \cdots \times \iota_{q_{\ell}}} & B{\rm SO}(q_1) \times \cdots \times B{\rm SO}(q_{\ell})   \ar[r] 
& \hspace{2mm}    B{\rm SO}(q)     
} 
\end{align}
For $q=2k$ and each $q_i = 2$, then on the bottom left side we get the product $\mC\mP^2 \times \cdots \times \mC\mP^2$
  of $k$ copies of $\mC\mP^2$. By the product formula \eqref{eq-multiplicative}, the pull-back of $p_1 \in H^4(B{\rm SO}(q) ; \mQ)$  is given by
  \begin{equation}\label{eq-p1}
p_1(\xi_1) \otimes 1 \otimes \cdots \otimes 1 + \cdots + 1\otimes   \otimes p_1(\xi_k) \in H^*(\mC\mP^2 ; \mQ) \otimes \cdots \otimes H^*(\mC\mP^2; \mQ) \ .
\end{equation}
In particular, the pull-back of $p_1^k$ is non-zero in $H^{4k}(\mC\mP^2 \times \cdots \times \mC\mP^2 ; \mQ)$, and so by the commutativity of the diagram \eqref{eq-commutingdiagram} we obtain that $\nu_q^*(p_1^k) \in H^{4k}(B\G_{2k}^+ ; \mQ)$ is non-zero.
  
  The diagram \eqref{eq-commutingdiagram} has another useful  consequence that we note. Given $q_1 > 1$ let $q_2 = q_3 = \cdots = q_{\ell+1} = 1$, then for $q= q_1 + \ell$, we obtain the commuting square
  \begin{align}\label{eq-commutingsquare}
\xymatrix{
  B\G_{q_1}^+   \ar[d]^{\nu_{q_1}}   \ar[r]
& \hspace{2mm} B\Gamma_{q_1+\ell}^+  \ar[d]^{\nu_{q_1+\ell}}  \\
  B{\rm SO}(q_1)    \ar[r] 
& \hspace{2mm}    B{\rm SO}(q_1+\ell)     
} 
\end{align}
so that injectivity results for the map $\nu_{q_1}^*$ imply injectivity for $\nu_{q_1 + \ell}^*$ for classes in the image of the restriction map $H^*(B{\rm SO}(q_1 + \ell) ; \mQ) \to H^*(B{\rm SO}(q_1) ; \mQ)$. In the following we abuse notation and identify the Pontrjagin classes in the latter space with those in the former space when appropriate.
   
 We  now give an extension of the   observation by Lawson. Recall that for each $1 < 2i< q$, the class $p_i \in H^{4i}(B{\rm SO}(q) ; \mQ)$ is supported by a spherical class in the image of the Hurewicz map $h \colon \pi_{4i}(B{\rm SO(q)}) \to H_{4i}(B{\rm SO}(q) ; \mQ)$. In particular, for $i \geq 2$ and $q=4i-2$, there exists a map $f_i \colon \mS^{4i} \to X(q) \to B\G_q^+$
   such that $f_i^*(\nu_q^* p_i) \ne 0$, and $f_i^*(\nu_q^* p_j) = 0$ for $0 < j < i$. 
  
  Given a collection of non-negative integers $\vec{n} = \{n_1 , \ldots , n_k\}$ set 
    $$p(\vec{n})   =  p_1^{n_1} \cdots p_k^{n_k} \in H^{q+ 2|\vec{n}|}(B{\rm SO}(q(\vec{n})) ; \mQ)$$
  where  $q(\vec{n}) = (4-2) n_1 + \cdots + (4k-2) n_k $ and   $|\vec{n}|  =  n_1 + \cdots + n_k$. Note that  $ q(\vec{n}) + 2|\vec{n}| \leq 2q(\vec{n})$.
 
 Here is the main result of this section.
 \begin{thm}\label{thm-independent}
For $q \geq 2$, define  
 \begin{equation}\label{eq-independent}
\fV(q) = \{ p(\vec{n})  \mid  k > 0, q(\vec{n}) \leq q\}  \subset H^*(B{\rm SO}(q) ; \mQ) \ .
\end{equation}
 Then   $\nu_q^* \fV(q) \subset H^*(B\G_q^+ ; \mQ)$ consists of linearly independent vectors.
 \end{thm}
 \proof
 
Suppose there exists a non-trivial relation between the elements of   $\nu_q^* \fV(q)$, then there is some degree $2 \leq  n \leq 2q$ and a non-trivial relation   of the form
\begin{equation}\label{eq-dependent}
0 = \sum_{i=1}^{\mu} \ a_i \cdot \nu_q^*(p(\vec{n}_i)) \in H^n(B\G_q^+ ; \mQ) \ , \ 0 \ne a_i \in \mQ \ ,
\end{equation}
where  each $\vec{n}_i = \{n_{i,1} , \ldots , n_{i,k_i}  \}$ with $q(\vec{n}_i) \leq q$ and the degree $q(\vec{n}_i)+2|\vec{n}_i| = n$ for all $ 1\leq i \leq \mu$.

We can assume the indices $\vec{n}_i$ are ordered so that $i < j$ implies $k_i \leq k_j$, and if $k_i = k_j$ then $n_{i,k_i} \leq n_{j,k_j}$. The proof is by induction on the lexicographical ordering on the vectors $\vec{n}_i$ appearing in the sum \eqref{eq-dependent}. Define the 
maximal index $k_{max} = \max \{k_i \mid 1 \leq i \leq \mu\}$.

Let $k_{max} =1$,  then by assumption that \eqref{eq-dependent} is a homogeneous sum,   it reduces to $a_1 \cdot \nu_q^* (p_1^{n_1}) = 0$ where $4 n_1 \leq 2q$. Recall that $\iota_2 \colon \mC\mP^2 \to B\G_2^+$ satisfies $0 \ne \iota_2^*(p_1) \in H^4(B\G_2^+ ; \mQ)$ so then by \eqref{eq-p1} we have $\times^{n_1} \iota_2 \colon \times^{n_1} \mC\mP \to B\G_q^+$ satisfies $(\times^{n_1} \iota_2)^*(\nu_q^* p_1^{n_k}) \ne 0$ hence $a_1=0$, a contradiction.

 Now we prove the inductive step. Assume that for $k_0 > 1$ there are no non-trivial relations of the form \eqref{eq-dependent} for $k_{max} < k_0$, and if $k_{max} = k_0$ then there are no non-trivial relations for indices with $n_{i,k_0} <   n_{\mu, k_0}$.

 Consider a relation \eqref{eq-dependent} with  $k_{max} = k_0$. 
   Let $\mu_0 \leq \mu$ be such that 
 $k_i = k_0$ for $\mu_0 \leq i \leq \mu$ and $k_i < k_0$ for $i < \mu_0$.  Divide the sum in \eqref{eq-dependent} into three terms:
 \begin{eqnarray*}
P_1 & = &   \sum_{i=1}^{\mu_0 -1}  \ a_i \cdot  p_i^{n_{i,1}} \cdots p_{k_i}^{n_{i,k_i}} \ , \ {\rm where} ~   k_i < k_0 \\
P_2 & = &  \sum_{i=\mu_0}^{\mu_1}  \ a_i \cdot  p_1^{n_{i,1}} \cdots p_{k_0}^{n_{i,k_0}}   \ , \ {\rm where} ~   n_{i,k_0} <   n_{\mu_1, k_0} \\
P_3 & = & \sum_{i=\mu_1}^{\mu}   \ a_i \cdot  p_1^{n_{i,1}} \cdots p_{k_0 -1}^{n_{i,k_0-1}}      \ , \ {\rm where} ~     {n_{i, k_0}} = {n_{\mu, k_0}} ~ {\rm for} ~ \mu_1 \leq i \leq \mu \  .
 \end{eqnarray*} 
\begin{equation}\label{eq-3sum}
0 = \nu_q^* \left( \sum_{i=1}^\mu \ a_i \cdot  p(\vec{n}_i)  \right)   =       \nu_q^*P_1 + \nu_q^*P_2 + \nu_q^* \left(P_3 \cdot p_{k_0}^{n_{\mu, k_0}} \right)
\end{equation}
We now construct a ``test cycle'' which we evaluate the expression \eqref{eq-3sum} against. 

Let $X = Y \times \{ \mS^{4k_0}\times \cdots \times \mS^{4k_0}\}$ with $n_{\mu, k_0}$ factors, and where 
\begin{equation}\label{eq-bigprod} 
Y = \{\mC\mP^2 \times \cdots \times \mC\mP^2\} \times \{\mS^8 \times \cdots \times \mS^8\} \times \cdots \times 
\{ \mS^{4(k_0-1)} \times \cdots \times \mS^{4(k_0-1)} \}
\end{equation}
and the number of factors  in \eqref{eq-bigprod} of $\mC\mP^2$ is $n_{\mu,1}$ and the number of factors of $\mS^{4k_j}$  is $n_{\mu,j}$ for $2 \leq j \leq k_0-1$. Then by the choice of the indices in $\vec{n}_{\mu}$ there exists a product of the map  
 $\iota_2 \colon \mC\mP^2 \to B\G_2^+$ and the maps 
$\{  f_i \colon \mS^{4i} \to B\G_{4i-2}^+  \mid 1 \mid 2 \leq i \leq \mu\}$ giving a product map 
$F_Y \times F_\mu \colon Y \times \mS^{4k_0} \to B\G_q$ such that  $F^*(\nu_q^* P_1) = 0$ as the pull-backs of the Pontrjagin classes in the factor $\mS^{4k_0}$ vanish, and   the restriction to the factor $Y$ vanishes for degree reasons.
Also,  $F^*(\nu_q^* P_2) = 0$ as the pull-backs of the Pontrjagin classes to    the restriction to the factor $Y$ vanishes for degree reasons. Then calculate 
$$F^*(\nu_q^* \left(P_3 \cdot p_{k_0}^{n_{\mu, k_0}} \right)) =  F_Y^*(\nu_q^* P_3) \ F_{\mu}^*(\nu_q^* (p_{k_0}^{n_{\mu, k_0}} ) ) \ .$$
The term $F_{\mu}^*(\nu_q^* (p_{k_0}^{n_{\mu, k_0}} ) ) \ne 0$ by the choice of the map $F_{\mu}$. Thus, when we evaluate the class in \eqref{eq-3sum} against the cycle defined by $F_Y \times F_{\mu}$ we obtain
$0 =  F_Y^*(\nu_q^* P_3)$. But by the induction hypothesis, this implies the coefficients in $P_3$ must vanish, contrary to assumption.
 \endproof

\section{Rigid Secondary classes} \label{sec-secondary}

 The survey by Lawson   \cite{Lawson1977} gives a nice overview  of the multiple approaches in the literature to defining the secondary classes of smooth foliations, by Bernshtein and Rozenfeld \cite{BernshteinRozenfeld1972,BernshteinRozenfeld1973},   Bott and Haefliger \cite{BottHaefliger1972}, and Kamber and Tondeur \cite{KT1974c,KT1975a}.  The author's survey \cite{Hurder2009} discusses the various approaches to showing that these invariants are non-trivial for various constructions of foliations. In this section, we give the construction of examples of foliated manifolds with non-trivial secondary classes, especially the rigid secondary classes mentioned in Section~\ref{sec-intro}, using the results on the non-triviality of the Pontrjagin classes in Section~\ref{sec-pontrjagin}.

The normal bundle $Q$ to a smooth-foliation $\F$, when restricted to a leaf $L_x$ of $\F$, has  a natural flat connection $\nabla^{L_x}$ defined by the leafwise parallel transport on $Q$ restricted to $L_x$. An    \emph{adapted connection}  $\nabla^{\F}$ on $Q \to M$ is a connection whose restrictions to leaves equals this natural flat connection. An adapted connection  need not be flat over $M$.  The connection data provided by $\nabla^{\F}$ can be thought of as a  ``linearization'' of the normal structure to $\F$ along the leaves.  Thus, $\nabla^{\F}$  captures aspects of the data provided by the Haefliger groupoid $\GF^r$ of $\F$;   it is a ``partial linearization'' of the   nonlinear data which defines the homotopy type of $B\GF$.

Denote by $I(\mathfrak{gl}(q, \mR))$ the graded ring of adjoint-invariant polynomials on the Lie algebra $\mathfrak{gl}(q, \mR)$ of the real general linear group $GL(q, \mR)$.  As a ring, $I(\mathfrak{gl}(q, \mR)) \cong \mR[c_1, c_2, \ldots , c_q]$  is a polynomial algebra on $q$ generators, where the $i^{th}$-Chern polynomial $c_i$ has  graded degree $2i$.  
 Associate to each generator $c_i$ the closed $2i$-form $c_i(\Omega({\nabla^{\F}})) \in \Omega_{deR}^{2i}(M)$, so that one obtains the Chern homomorphism 
$\Delta^{\F} \colon \mR[c_1, c_2, \ldots , c_q] \to \Omega_{deR}^{ev}(M)$.     
       The   Bott Vanishing Theorem holds at the level of the differential forms representing the characteristic classes of the normal bundle to a $C^2$-foliation, which is the basis for the construction  of the secondary characteristic classes for $C^2$-foliations.

 \begin{thm}[Strong Bott Vanishing \cite{BottHaefliger1972}] \label{SBVT} Let $\F$ be a codimension-$q$, $C^2$-foliation, and $\nabla^{\F}$ an adapted connected on the normal bundle $Q$.   Then for any polynomial $c_J \in \mR[c_1, c_2, \ldots , c_q]$ of graded degree $deg(c_J) > 2q$, the Chern form $c_J(\Omega({\nabla^{\F}}))$ is identically zero. 
 \end{thm}
 
 Thus, there is an induced map of DGA's (differential graded algebras):  
 \begin{equation}\label{eq-SBVT}
\Delta^{\F} \colon   \mR[c_1, c_2, \ldots , c_q]_{2q} \longrightarrow \Omega_{deR}^{ev}(M)  \ ,
\end{equation}
where $\mR[c_1, c_2, \ldots , c_q]_{2q}$ denotes the polynomial ring truncated in degrees greater than $2q$.

Now assume that the normal bundle to the foliation $\F$ on $M$ is trivial. The choice of a framing of $Q$, denoted by $s$, induces an isomorphism $Q \cong M \times {\mathbb R}^q$. Let $\nabla^s$ denoted the connection for the trivial bundle, then the curvature forms of $\nabla^s$ vanish identically. It follows that the transgression form $y_i =Tc_i(\nabla^{\F}, \nabla^s) \in \Omega_{deR}^{2i-1}(M)$   of $c_i$   satisfies the equation 
$dy_i = c_i(\Omega({\nabla^{\F}}))$.
Consider the DGA   
\begin{equation}\label{eq-defWq}
W_q  = \Lambda(y_1, y_2, \ldots, y_{q}) \otimes \mR[c_1, c_2, \ldots , c_q]_{2q} \ ,
\end{equation}
where the differential is defined by $d(y_i\otimes 1) = 1 \otimes c_i$.

The data $(\F, s,\nabla^{\F})$  determine a map $\Delta^{\F,s} \colon W_q \to \Omega^*(M)$ of DGA's.
The induced map in cohomology, $\ds \Delta^{\F,s}  \colon H^*(W_q) \to H^*(M)$,   depends only on the homotopy class of the framing $s$ and the framed concordance class of $\F$, and moreover, this construction is functorial.

\begin{thm}    There is a    well-defined   universal characteristic map
\begin{equation}
\Delta   \colon    H^*(W_q)    \to    H^*(B\overline{\Gamma}_q ; \mR) 
\end{equation}
Given a   codimension-$q$ foliation $\F$ with framing $s$, the classifying map $\ds \Delta^{\F,s} \colon H^*(W_q) \to H^*(M; \mR)$ factors through the universal map:

\begin{picture}(100,45)\label{comm_diag}
\put(170,40){$H^*(B\overline{\Gamma}_q ; \mR) $}
\put(190,30){\vector(0,-1){17}}
\put(130,13){\vector(3,2){30}}
\put(200,20){$h_{\F}^*$}
\put(130,25){$\Delta$}
\put(150,8){$\Delta^{\F, s}$}
\put(100,0){$H^*(W_q)$}
\put(145,2){$\vector(1,0){22}$}
\put(176,0){$H^*(M; \mR)$}
\end{picture}
\end{thm}

We recall the  \emph{Vey basis} for $H^{*}(W_q)$, as first described  in \cite{Godbillon1974}:
\begin{prop}\label{prop-veybasisWq}
The following set of monomials in $W_q$ form a basis for $H^{*}(W_q)$:
\begin{eqnarray} \label{eq-VeybasisWq}
y_I c_J  & \ {\rm such \  that} \ &  I = (i_1, \ldots , i_s) \ {\rm with} \ 1 \leq i_1 < \cdots < i_s \leq q      \\
  &   &   J = (j_1, \ldots , j_{\ell}) \ {\rm with} \  j_1 \leq \cdots \leq j_{\ell} \ , \  j_1 + \cdots j_{\ell} \leq q \nonumber\\
    &   &    i_1 + j_1 + \cdots j_{\ell} \geq q+1 \ , \  i_1 \leq    j_1 \ .  \nonumber
\end{eqnarray}
\end{prop}

   We also briefly recall the construction of the secondary classes for a   smooth foliation $\F$ of a manifold $M$, with no assumption that its normal bundle is trivial. In place of the flat connection $\nabla^s$ associated to a framing of $Q$, let $\nabla^g$ be the connection on $Q$ associated to a Riemannian metric on $Q$. Then the closed forms $c_{2i}(\Omega(\nabla^g)) \in \Omega^{4i}_{deR}(M)$ need not vanish, and in fact their cohomology classes define the Pontrjagin forms for the normal bundle.  On the other hand,  the forms $c_{2i+1}(\Omega(\nabla^g)) \in \Omega^{4i+2}_{deR}(M)$ vanish by the skew-symmetry of the curvature matrix $\Omega(\nabla^g)$.
  Thus, we can repeat the above construction of a characteristic map for the subcomplex $W_q$ of $W_q$ defined by
    \begin{equation}
WO_q  = \Lambda(y_1, y_3, \ldots, y_{q'}) \otimes \mR[c_1, c_2, \ldots , c_q]_{2q} \ , 
\end{equation}
where $q' \leq q$ is the largest odd integer with $q' \leq q$.

  The data $(\F, \nabla^{g})$  determine a map    $ \Delta^{\F} \colon WO_q \to \Omega^*(M)$ of DGA's.
The induced map in cohomology, $\ds \Delta^{\F}  \colon H^*(WO_q) \to H^*(M)$,   depends only on the   framed concordance class of $\F$,  but not on the choice of the Riemannian metric on $Q$. It follows that there is a   well-defined   universal characteristic map
$\Delta   \colon    H^*(WO_q)    \to    H^*(B\Gamma_q ; \mR)$.

The   \emph{secondary characteristic classes} of   foliations are  defined as  those in the image of the maps  $\Delta$ for degree at least $2q+1$, in either the framed or unframed cases.

  \begin{defn}\label{def-specialcases} 
  The rigid secondary classes for foliations with framed normal bundles are those in the image of the restriction map $H^*(W_{q+1}) \to  H^*(W_q)$.    In terms of the Vey basis, the class  $h_I \wedge c_J$ is said to be \emph{rigid} if $i_1 + j_1 + \cdots j_{\ell} \geq q+2$.
  \end{defn}
 
We next give the   construction of examples with non-trivial rigid secondary classes.
Let $X$ be a compact manifold and $f_Q \colon X \to B{\rm SO}(q)$ classify a bundle $Q \to X$.
Suppose there exists a lifting of $f_Q$ to $\nu_f \colon X \to B\G_q^+$, then $\nu_Q$ determines a $\G_q^+$-structure on $X$ with normal bundle   $Q$.
Let $M = {\rm Fr}(Q) \to X$ denote the ${\rm SO}(q)$-bundle of orthogonal frames for $Q$, so that we obtain a fibration
${\rm SO}(q) \to M \to X$. 
The lift of the $\G_q^+$-structure on $X$ to $M$ then has a canonical framing $s$, so is classified by a map $\nu_f^s \colon M \to B\oG_q$. Then there is a commutative diagram, where each vertical sequence is a principal ${\rm SO}(q)$-fibration:
  \begin{align}\label{eq-pullback}
\xymatrix{
  {\rm SO}(q)   \ar[d]^{\iota}   \ar[r]^{id}  & \hspace{2mm}     {\rm SO}(q)  \ar[d]^{\iota}  \ar[r]^{id} 
& \hspace{2mm} {\rm SO}(q)  \ar[d]^{\iota}  \\
 M   \ar[d]^{\pi}   \ar[r]^{f_Q^s} & \hspace{2mm} B\oG_q  \ar[d]   \ar[r]^{\widetilde{\nu}}   & \hspace{2mm}    E{\rm SO}(q)   \ar[d]^{\pi}  \\
  X    \ar[r]^{f_Q}  & \hspace{2mm}    B\G_q^+  \ar[r]^{\nu}   & \hspace{2mm}    B{\rm SO}(q) 
} .
\end{align}
For $q=2k+1$, the Chern-Weil model for the cohomology of the acyclic space $E{\rm SO}(q)$ is given by    
\begin{equation}\label{eq-CWSO}
\Lambda(\widehat{p}_1, \ldots , \widehat{p}_k)  \otimes {\mathbb R}[p_1, \ldots , p_k]    \quad , \quad d(\widehat{p}_i \otimes 1) = 1 \otimes p_i \ ,
\end{equation}
  For $q$ even, one modifies this complex by adding the Euler class $e_q$ and its transgression class $\widehat{e}_q$ with $d(\widehat{e}_q \otimes 1) = 1 \otimes e_q$. By the naturality of the Chern-Weil construction, the pull-back of the DGA model in \eqref{eq-CWSO} to the middle fibration is the DGA subalgebra of the DGA model in \eqref{eq-defWq},
 \begin{equation}\label{eq-subalgebra}
  \Lambda(y_2, \ldots, y_{2k}) \otimes \mR[c_2, \ldots , c_{2k}]_{2q} \subset \Lambda(y_1, y_2, \ldots, y_{q}) \otimes \mR[c_1, c_2, \ldots , c_q]_{2q} = W_q\ .
\end{equation}
 Finally, the pull-back of the forms $\widehat{p}_i$ to the left hand fibration gives the forms $\Delta^{\F,s}(y_{2i}) \in \Omega^{4i-1}(M)$, and the pul-back of the forms $p_i$ gives $\nu_f^*(p_i) \in \Omega^{4i}(M)$.
 We use this to calculate the image of the characteristic map $\ds \Delta^{\F,s} \colon H^*(W_q) \to H^*(M; \mR)$ for various choices of spaces $X$ and lifted maps $f_Q^s$.
 
The first construction of examples with non-trivial rigid secondary classes  is for even codimension $q=2k \geq 4$.
Let $\mC\mP^{\infty}$ denote the infinite projective space, which serves as a model for $B{\rm SO}(2)$, and let  $\xi \to \mC\mP^{\infty}$ be the canonical $\mR^2$-bundle.  The Euler class $e(\xi) \in H^2(B{\rm SO}(2) ; \mZ)$ is  a generator of the ring $H^*(B{\rm SO}(2) ; \mZ)$, and the first Pontrjagin class $p_1(\xi) = e(\xi)^2$   is a generator of $H^4(\mC\mP^2 ; \mZ)$. 
The 4-skeleton of a CW decomposition of $\mC\mP^{\infty}$ is the complex projective space   $\mC\mP^2$, and we also denote by $\xi \to  \mC\mP^2$   the restriction of the canonical bundle.

 Let $X_k = \mC\mP^2 \times \cdots \times \mC\mP^2$ the $k$-fold product, and let $f_Q \colon X_k \to B{\rm SO}(2k)$ be the product of the inclusions followed by the map
 $B{\rm SO}(2) \times \cdots \times B{\rm SO}(2) \subset B{\rm SO}(2k)$. That is, $f_Q$ classifies the direct sum 
    $ \Xi_k  = \xi_1 \oplus \cdots \oplus \xi_k$, where $\xi_{i} \to X_k$ denotes the canonical bundle over the $i$-th factor in $X_k$.
  Then the product formula for the rational Pontrjagin classes yields
\begin{eqnarray}
p_{\ell}(\Xi_k) & = &  \sum  \  p_{i_1}(\xi_1) \cdots p_{i_n}(\xi_k)  \nonumber  \\
  & = &  \sum_{i_1 < \cdots <  i_{\ell}} \  p_{1}(\xi_{i_1}) \cdots p_{1}(\xi_{i_\ell})  \in H^{4 \ell}(X_k ; \mQ)  \label{eq-product}
\end{eqnarray}
where we use that all products of Pontrjagin classes in $H^*(X_k ; \mQ)$ vanish. In particular, \eqref{eq-product} implies that for $\ell > 1$ the class $p_{\ell}(\Xi_k)$ is a rational multiple of $p_1(\Xi_k)^{\ell}$.

The fiber $B\oG_2$ is $3$-connected, so there is a lifting $\widetilde{\iota} \colon \mC\mP^2 \to B\G_2^+$ of the inclusion map $\mC\mP^2  \subset B{\rm SO}(2)$. Let  $\widetilde{f}_Q \colon X_k \to B\G_{2k}^+$ be lifting of the map $f_Q$ obtained by taking products of the map $\widetilde{\iota}$.

Form the pull-back  bundle  ${\rm SO}(2k) \to M_k \to X_k$ as in the diagram \eqref{eq-pullback}, then the canonical framing of the pull-back of $Q$ to $M_k$ determines classifying map $\widetilde{f}_Q^s \colon M_k \to B\oG_{2k}$.
We calculate the image of the classifying map $\Delta^{\widetilde{f}_Q^s} \colon H^*(W_{2k}) \to H^*(M_k; \mR)$ using the functoriality of Chern-Weil construction.
Let $\{x_1, \ldots , x_k\} \subset H^2(X_k, \mR)$ be the generators corresponding to the factors of $X_k$, then in notation as above we have $p_1(\xi_i) = x_i^2$, so $f_Q^*(p_1) = x_1^2 + \cdots + x_k^2$. 

Consider then the $E_2$-term of the spectral sequence for the fibration ${\rm SO}(2k) \to M_k \to X_k$,
\begin{equation}
E_2^{r,s} = H^r(X_k ;  H^s({\rm SO}(2k) ; \mR)) \cong   H^r(X_k ; \mR) \otimes H^s({\rm SO}(2k) ; \mR) \ .
\end{equation}

 By diagram \eqref{eq-pullback} and the functoriality of the construction of the spectral sequence,     the map on $E_2$-terms induced by    $\Delta^{\widetilde{f}_Q^s}$ is the identity on the fiber cohomology, and induced by $f_Q^*$ on the base cohomology.
 That is, $u_i = \Delta^{\widetilde{f}_Q^s}(y_{2i})$ is the primitive generator   $H^{4i-1}({\rm SO}(2k) ; \mR)$, and so the $E_2$-differential   satisfies $d_2(1 \otimes u_i) = f_Q^*(p_i) \otimes 1$. A straightforward calculation using \eqref{eq-product} then yields:
 
 \begin{prop}\label{prop-2k-indepedent}
 For $I=(2=2i_1 < 2i_2 < \cdots 2i_{\ell})$ the cohomology classes
  \begin{equation} \label{eq-2k-indepedent}
\Delta^{\widetilde{f}_Q^s}(y_I \wedge c_2^k) = f_Q^*(p_1)^k \otimes u_{i_1} \wedge \cdots u_{i_{\ell}} \in H^*(M_k; \mR)
\end{equation}
are linearly independent.
 \end{prop}
  
The second construction of examples with non-trivial rigid secondary classes is for even codimension $q=4k-2 \geq 6$.   These classes are independent of the classes in Proposition~\ref{prop-2k-indepedent}. 
 
It was noted in the previous section that for $1 \leq i < q/2$, the Pontrjagin class   $p_i$ pairs non-trivially with a class in the image of the Hurewicz map $h \colon \pi_{4i}(B{\rm SO}(q)) \to H_{4i}(B{\rm SO}(q) ; \mZ)$. Let $g_Q \colon \mS^{4k} \to B{\rm SO}(q)$ be such that $g_Q^*(p_k) \ne 0$. Note that $g_Q^*(p_i) = 0$ for $i \ne k$, and in particular $g_Q^*(p_1) = 0$.

The fiber $B\oG_q$ is $(q+1)$-connected, so there is a lifting $\widetilde{g}_Q \colon \mS^{4k} \to B\G_q^+$ of the  map $g_Q \colon \mS^{4k}  \to B{\rm SO}(q)$. 

Let ${\rm SO}(q) \to N_k \to \mS^{4k}$ be the pull-back of the principle ${\rm SO}(q)$-bundle over $B{\rm SO}(q)$, then $\widetilde{g}_Q$ determines a   $\G_q^+$-structure on $N_k$ with a canonical framing of its normal bundle, hence  is classified by a map $\widetilde{g}_Q^s \colon N_k \to B\oG_q$ that is a lift of $\widetilde{g}_Q$.

 Again,   calculate the image of the classifying map $\Delta^{\widetilde{g}_Q^s} \colon H^*(W_{q}) \to H^*(N_k; \mR)$ using the functoriality of Chern-Weil construction.
Let $\chi_{4k} \in H^{4k}({\rm S}^{4k} ; \mZ)$ be a generator, then   $\widetilde{g}_Q^*(p_k) = c_k \cdot \chi_{4k}$ for some constant $c_k \ne 0$. 
For the $E_2$-term of the spectral sequence for the fibration ${\rm SO}(q) \to N_k \to {\rm S}^{4k}$,
\begin{equation}
E_2^{r,s} = H^r({\rm S}^{4k} ;  H^s({\rm SO}(q) ; \mR)) \cong   H^r( {\rm S}^{4k}; \mR) \otimes H^s({\rm SO}(q) ; \mR) \ ,
\end{equation}

 By diagram \eqref{eq-pullback} and the functoriality of the construction of the spectral sequence,     the map on $E_2$-terms induced by    $\Delta^{\widetilde{g}_Q^s}$ is the identity on the fiber cohomology, and induced by $g_Q^*$ on the base cohomology.
 That is, $u_i = \Delta^{\widetilde{g}_Q^s}(y_{2i})$ is the primitive generator   $H^{4i-1}({\rm SO}(q) ; \mR)$, and so the $E_2$-differential  satisfies $d_2(1 \otimes u_i) = g_Q^*(p_i) \otimes 1$. A straightforward calculation using \eqref{eq-product} then yields:
 
 \begin{prop}\label{prop-4k-2-independent}
 For $q=4k-2 \geq 6$ and $I=(2k=2i_1 < 2i_2 < \cdots 2i_{\ell})$ the   classes
  \begin{equation} \label{eq-2k-2-indepedent}
\Delta^{\widetilde{g}_Q^s}(y_I \wedge c_{2k}) = g_Q^*(p_k)  \otimes u_{i_1} \wedge \cdots u_{i_{\ell}} \in H^*(N_k; \mR)
\end{equation}
are linearly independent.
 \end{prop}
It follows that the classes $\Delta(y_I \wedge c_{2k}) \in H^*(B\oG_q  ; \mR)$ are linearly  independent. Also, $g_Q^*(p_1) = 0$ implies that 
$\Delta^{\widetilde{g}_Q^s}(y_2 \wedge c_2^{2k-1}) =0$, so these classes are also independent of $\Delta(y_2 \wedge c_2^{2k-1})$.

The third construction   is based on the    \emph{permanence principle}, as it was called by Lazarov in \cite{Lazarov1979} and extended by the author in \cite{Hurder1981b}. In this construction, we use the ${\rm SO}(q)$-action fiberwise on $B\oG_q$ to  extend the non-triviality of one secondary class, to a family of non-trivial classes for $B\oG_q$. We recall a simple  version of this construction.
  
  Let  $f_{I,J} \colon \mS^n \to B\oG_q$ be a map such that $f_{I,J}^*(y_I c_J) \in H^n(\mS^n; \mR)$ is non-vanishing, where $y_I c_J$ is an element of the Vey basis for $H^*(W_q)$. 
  
  The normal bundle $f_{I,J}^! Q$ to the induced $\oG_q$-structure on $\mS^n$ has a framing $s_{I,J}$, which induces a framing $s_P$ on    the pull-back bundle $Q_P \to P$, for the product space $P =  {\rm SO}(q) \times \mS^n$. Let $s_P'$ denote the framing of $Q_P$ obtained by twisting $s_P$ over the factor ${\rm SO}(q)$ using the canonical action on frames. Then by \cite[Theorem~3.3]{Hurder1981b} we have:

   \begin{prop}\label{prop-permanance}
 For $I= (i_1 < i_2 < \ldots < i_{\ell})$, $K=(2k_1 < 2k_2 < \cdots < 2k_{\mu})$ with $i_{\ell} < 2k_1$, the cohomology classes
  \begin{equation} \label{eq-2k-permanance}
\Delta^{f_{I,J}^{s_P'}}(y_I \wedge y_K \wedge c_J)   \in H^*(P; \mR)
\end{equation}
are linearly independent.
 \end{prop}

\section{Spherical classes}\label{sec-spherical}

Given a foliation $\F$ of codimension $q$ on a paracompact manifold $M$, there is a classifying map $\nu_{\F} \colon M \to B\G_q$ whose homotopy class is uniquely defined by $\F$. Conversely, a  map $f \colon M \to B\G_q$ determines a     $\G_q$-structure on $M$, which is a foliated microbundle over $M$ (see \cite{Haefliger1970,Lawson1977}).   This motivates the interest in the set $[M, B\G_q]$ of homotopy classes of maps from $M$ to $B\G_q$, and similarly for the set $[M, B\oG_q]$ which classifies   $\G_q$-structures on $M$ with framed normal bundles. 
The following result from \cite[proof  of Theorem 4.5]{Hurder1981b}  is useful for the study of $[M, B\oG_q]$.
\begin{prop}\label{prop-BP}
Let $M$ be a compact oriented $n$-manifold. Then there is an exact sequence of sets
\begin{equation}\label{eq-BP}
 [\Sigma M' , B\oG_q] \stackrel{\beta'}{\longrightarrow}  \pi_n(B\oG_q) \stackrel{\rho'}{\longrightarrow} [M, B\oG_q] \ ,
\end{equation}
where the image of $\beta'$ is contained in the torsion subgroup of $\pi_n(B\oG_q)$.
\end{prop}
Here, the space $M'$ is the (n-1)-skeleton for a CW decomposition of $M$, $\rho \colon M \to M/M' \cong \mS^n$ is the map collapsing $M'$ to a point, $\rho'$ is the induced map, and $\beta'$ is the map induced by the connecting map $\beta \colon \mS^n \to \Sigma M'$ for the Barratt-Puppe sequence of the attaching map $\alpha \colon \mS^{n-1} \to M'$.
See the details of the proof of \cite[Theorem 4.5]{Hurder1981b} for further details.
Observe that  a real valued invariant  for $\pi_n(B\oG_q)$ must vanish on the image of $\beta'$,  hence any element of  $\pi_n(B\oG_q)$ which is non-trivial for the invariant determines a non-trivial $\oG_q$-structure  on $M$.

The Rational   Hurewicz Theorem is the basic result for the study of spherical cohomology:  
\begin{thm}[see \cite{KlausKreck2004}] \label{thm-RHIT}
Let $Y$ be a simply connected space with $\pi_k(Y) \otimes \mQ = 0$ for $0 \leq k \leq r$. Then the Hurewicz map induces an isomorphism  
\begin{equation}
h \otimes \mQ \colon \pi_{k}(Y) \otimes \mQ  \longrightarrow   H_{k}(Y ; \mQ)
\end{equation}
for $1 \leq k \leq 2r$, and a surjection for $k=2r+1$.
\end{thm}
 
 We next use Theorem~\ref{thm-RHIT} to show that $H^n_s(B\oG_q ; \mR)$ is non-trivial for many degrees $n >2q$.
 
Let $q=2k \geq 4$ and $M_k$ be the space constructed in the proof of Proposition~\ref{prop-2k-indepedent}, and $\widetilde{f}_Q^s \colon M_k \to B\oG_{2k}$ the classifying map. Then the rigid class $\Delta^{\widetilde{f}_Q^s}(y_2 \wedge c_2^k)   \in H^{2+3}(M_k; \mR)$ is non-trivial. Now let $\overline{M}_k$ denote the quotient space of $M_k$ by its (q+1)-skeleton, for some CW decomposition of $M_k$.   As $B\oG_q$ is (q+1)-connected, the map $\widetilde{f}_Q^s$ descends to a map $\overline{f}_Q^s  \colon \overline{M}_k \to B\oG_q$ such that $0 \ne \Delta^{\overline{f}_Q^s}(y_2 \wedge c_2^k)   \in H^{2q+3}(\overline{M}_k; \mR)$.
The map $h \colon \pi_{4k+3}(\overline{M}_k) \otimes \mQ  \longrightarrow   H_{4k+3}(\overline{M}_k ; \mQ)$ is onto by Theorem~\ref{thm-RHIT}, so there exists a map $\xi_k \colon \mS^{4k+3} \to \overline{M}_k$ such that 
$\langle \Delta^{\overline{f}_Q^s}(y_2 \wedge c_2^k) , h[\xi_k] \rangle \ne 0$. Thus the class $ \Delta(y_2 \wedge c_2^k) \in H^{q+3}(B\oG_q ; \mR)$ pairs non-trivially with the spherical cycle $h [\overline{f}_Q^s  \circ \xi_k ] \in H_{q+3}(B\oG_q ; \mZ)$.

Next, let $q = 4k -2 \geq 6$ and let $N_k$ be the  space constructed in the proof of Proposition~\ref{prop-4k-2-independent}, and $\widetilde{g}_Q^s \colon N_k \to B\oG_{q}$ the classifying map. Then the rigid class $\Delta^{\widetilde{g}_Q^s}(y_{2k} \wedge c_{2k})   \in H^{2q+3}(M_k; \mR)$ is non-trivial. Now proceed as for the case $q=2k$ above.
Let $\overline{N}_k$ denote the quotient space of $N_k$ by its (q+1)-skeleton, for some CW decomposition of $N_k$.   As $B\oG_q$ is (q+1)-connected, the map $\widetilde{g}_Q^s$ descends to a map $\overline{g}_Q^s  \colon \overline{N}_k \to B\oG_q$ such that $0 \ne \Delta^{\overline{g}_Q^s}(y_{2k} \wedge c_{2k})   \in H^{2q+3}(\overline{N}_k; \mR)$.
The map $h \colon \pi_{4k+3}(\overline{N}_k) \otimes \mQ  \longrightarrow   H_{4k+3}(\overline{N}_k ; \mQ)$ is onto by Theorem~\ref{thm-RHIT}, so there exists a map $\eta_k \colon \mS^{4k+3} \to \overline{N}_k$ such that 
$\langle \Delta^{\overline{g}_Q^s}(y_{2k} \wedge c_{2k}) , h[\eta_k] \rangle \ne 0$. 
Thus the class $ \Delta(y_{2k} \wedge c_{2k}) \in H^{q+3}(B\oG_q ; \mR)$ pairs non-trivially with the spherical cycle $h[\overline{g}_Q^s  \circ \eta_k ] \in H_{q+3}(B\oG_q ; \mZ)$.

Finally, we give a form of the permanence principle for spherically supported secondary classes.  Suppose there is map $f \colon \mS^{n} \to B\oG_q$ such that 
$f^*(\Delta(y_I c_J)) \ne 0$  for   $y_I c_J$ an element of  the Vey basis of $H^*(W_q)$, where 
$I=(i_1 < \dots < i_{\ell})$. 
Apply Theorem~\ref{prop-permanance} for   $f \colon \mS^n \to B\oG_q$ to obtain  the ${\rm SO}(q)$-twisted classifying map $\tau \colon P = {\rm SO}(q) \times \mS^n \to B\oG_q$.
Then for    $r$ with $i_{\ell} < 2r \leq q+2$,  
$$0 \ne \tau^*(\Delta(y_I \wedge y_{2r} \wedge c_J)) \in H^*(P ; \mR)  \ ,$$
by Proposition~\ref{prop-permanance} in \cite{Hurder1981b}. Moreover, the proof shows that  for the isomorphism  $H^*(P ; \mR) \cong H^*({\rm SO}(q) ; \mR) \otimes H^*(\mS^n ; \mR)$, we have
$\tau^*(\Delta(y_I \wedge y_{2r} \wedge c_J)) = \pm Tp_{r} \otimes f^*(\Delta(y_I c_J))$. 
 Choose a map $g_{r} \colon \mS^{4r-1} \to {\rm SO}(q)$ such that $g_{r}^*(Tp_{r}) \ne 0$.
Then for the composition
$$\mS^{4r-1} \times \mS^n \stackrel{g_r \times id}{\longrightarrow} {\rm SO}(q) \times \mS^n \stackrel{\tau}{\longrightarrow} B\oG_q$$
 we have $\{(g_r \times id) \circ \tau\}^*(\Delta(y_I \wedge y_{2r} \wedge c_J)) \ne 0$. The assumption that $2r \leq q+2$ implies that the image of $(g_r \times id) \circ \rho(\mS^{4r-1} \times pt) \subset B\oG_q$ is contractible. Let 
 $$Y= (\mS^{4r-1} \times \mS^n)/(\mS^{4r-1} \times pt) \cong \mS^{n+4r-1} \vee \mS^n \ ,$$
 then $(g_r \times id) \circ \tau$ descends to a map $\xi \colon Y \to B\oG_q$ such that $\xi^*(\Delta(y_I \wedge y_{2r} \wedge c_J)) \ne 0$. As this class must be supported on the first term $\mS^{n+4r-1}$, we have that  
 $\Delta(y_I \wedge y_{2r} \wedge c_J) \in H_s^{n+4r-1}(B\oG_q ; \mR)$.
 
 Introduce the following  collections of secondary classes:
 \begin{eqnarray}
 q=4k  \geq 4, & \ocR_q^s   = & \{y_2 \wedge y_K \wedge c_2^k \mid K = (2k_1 < \cdots < 2k_{\ell}) , 1< k_1, k_{\ell} \leq  k\} \label{eq-4k-rigid}\\
 q=4k-2 \geq 6, & \ocR_q^s   =& \{y_{2k} \wedge c_{2k},  y_2 \wedge y_K \wedge c_2^k \mid K = (2k_1 < \cdots < 2k_{\ell}) ,  1< k_1, k_{\ell} \leq k \} \ . \label{eq-4k-2-rigid}
 \end{eqnarray}
  Apply the above procedure iteratively to the class  $[\overline{f}_Q^s  \circ \xi_k ] \in \pi_{q+3}(B\oG_q)$   constructed above to obtain the proof of Theorem~\ref{thm-rigidspherical}.

 Proposition~\ref{prop-BP} can also can be applied for the dual homotopy invariants of foliations defined in \cite{Hurder1981a}, to obtain ``rigid homotopy classes'' in $\pi_*(B\oG_q)$ for degrees tending to infinity, as was shown in \cite{Hurder1985a}.

\section{Applications}\label{sec-applications}

In this section, we apply the conclusions of Theorem~\ref{thm-rigidspherical} to construct solutions of Problem~\ref{prob-homotopy}. The results  of Theorem~\ref{thm-rigidspherical} are also applied  in \cite{Hurder2024b} to show that     the homotopy groups $\pi_*(B\G_q)$  contain divisible subgroups,    and applied  in \cite{Hurder2023} to construct  group actions with non-trivial secondary invariants that are invariant under deformations of the group actions.

Let $M$ be a compact $n$-manifold without boundary. Suppose that the tangent bundle admits a decomposition $TM = F \oplus Q$\, classified by a map $(\tau_F, \nu_Q) \colon M \to B{\rm O}(n-q) \times B{\rm O}(q)$. Given a lifting $\widetilde{\nu}_Q \colon M \to B\G_q$ of the map $\nu_Q$, Thurston proved in \cite{Thurston1974b} that there is a codimension $q$ foliation $\F$ of $M$ whose normal bundle is homotopic with $Q$, and its classifying map is homotopic with $\widetilde{\nu}_Q$.  We call this the ``Thurston Realization Theorem''.  An alternate proof of it   was given by  Mi\v{s}a\v{c}ev  and   \`Elia\v{s}berg  in 
\cite{MisacevEliasberg1977,MisacevEliasberg1997}.  Both approaches to the proof are ``non-constructive'', in that the foliation is  shown to exist, but there is no information obtained about the geometric properties of the foliation.

We apply the Thurston Realization Theorem  for the case when $Q$ is a trivial bundle, so the classifying map restricts to a map $\widetilde{\nu}_Q \colon M \to B\oG_q$.
There are two cases to consider, when $q=2k \geq 4$ and $q=4k-2 \geq 6$. We illustrate the construction in the lowest codimension for each case.

Let $q = 4= 2k$, then $\ocR^s_4 = \{y_2 c_2^2\}$ contains one class of degree $n = 2q+3 = 11$.
Let $M$ be a compact orientable $11$-manifold without boundary. We assume in addition that $TM$ contains a trivial rank-4 bundle $\e_4 \subset TM$. For example, let  $M = M_0 \times \mT^4$ where $\mT^4$ is the 4-torus, and $M_0$ is any compact oriented $7$-manifold. Another interesting example is $M = \mS^{10} \times \mS^1$, which is a parallelizable manifold and in particular contains a trivial rank-4 subbundle.
 
 Let $M_2$ be the space constructed in the proof of Proposition~\ref{prop-2k-indepedent}, and $\widetilde{f}_Q^s \colon M_2 \to B\oG_{4}$ the classifying map. Let   $\overline{f}_Q^s  \colon \overline{M}_2 \to B\oG_4$ such that $0 \ne \Delta^{\overline{f}_Q^s}(y_2 \wedge c_2^2)   \in H^{11}(\overline{M}_2; \mR)$ be the map constructed above, 
 and  $\xi_2 \colon \mS^{11} \to \overline{M}_2$ the map such that   the class $ \Delta(y_2 \wedge c_2^2) \in H^{11}(B\oG_4 ; \mR)$ pairs non-trivially with the spherical cycle $h  [\overline{f}_Q^s  \circ \xi_2 ] \in H_{11}(B\oG_4 ; \mZ)$.  
 
Let $M$ be the given manifold, and $\rho' \colon \pi_{11}(B\oG_4) \to [M, B\oG_4]$ the map constructed in Proposition~\ref{prop-BP}, which was defined using the map $\beta \colon M \to \mS^{11}$ obtained by collapsing onto the top cell of $M$. For  each integer $\ell \in \mZ$, let $f_{\ell} \colon \mS^{11} \to B\oG_4$ be a map representing the homotopy class $\ell \cdot [\overline{f}_Q^s  \circ \xi_2 ] \in \pi_{11}(B\oG_4)$. 
Then the values
\begin{equation}
\langle (f_{\ell} \circ \beta)^*(\Delta(y_2 \wedge c_2^2)) , [M] \rangle = \ell \cdot \langle \Delta(y_2 \wedge c_2^2) , h  [\overline{f}_Q^s  \circ \xi_2 ] \rangle \ne 0
\end{equation}
  are distinct for all $\ell$, and so for $\ell \ne \ell'$ the classifying maps $f_{\ell} \circ \beta \colon M \to B\oG_4$ are pairwise non-homotopic. Let $\F_{\ell}$ be a codimension 4 foliation on $M$ constructed using the Thurston Realization Theorem, then 
 $\F_{\ell}$ and $\F_{\ell'}$ are not homotopic foliations for $\ell \ne \ell'$. On the other hand, the foliations $\F_{\ell}$ on the fixed manifold $M$ have homotopic normal bundles, and so homotopic tangent bundles.
 
 Next consider the case when $q=6 = 4k-2$, then $\ocR^s_6 = \{y_2 c_2^3, y_4c_4, y_2 y_4 c_2^3\}$ where the first two classes have degree $n = 2q+3 = 15$, and the third has degree $n=2q+10 = 22$.  First, let $M$ be a compact orientable $15$-manifold without boundary, such  that $TM$ contains a trivial rank-6 bundle. 
 In this particular case, we can take $M = \mS^{15}$ which by the results of  James \cite{James1957} contains a trivial subbundle of rank $8$, hence a trivial rank 6 subbundle. On the other hand, the celebrated work of  Adams \cite{Adams1962} shows that this is the unique case for which our results apply. Alternately,  let $M = \mS^{14} \times \mS^1$ for example, which is a parallelizable manifold. There are many further possibilities.
 
By the results of Section~\ref{sec-spherical},  there exists   $f = \overline{f}_Q^s  \circ \xi_2 \colon \mS^{15} \to B\oG_6$ and 
$g = \overline{g}_Q^s  \circ \eta_2 \colon \mS^{15} \to B\oG_6$
with 
 \begin{equation}
\langle f^*(y_2 c_2^3) , [M] \rangle \ne 0 ; \langle g^*(y_2 c_2^3) , [M] \rangle = 0 ; \langle g^*(y_4c_4) , [M] \rangle \ne 0 \ ,
\end{equation}
so in particular their homotopy classes $\{[f], [g]\} \subset \pi_{15}(B\oG_6) \otimes \mR$ are linearly independent.
 
 For each pair of integers $(\ell_1, \ell_2) \in \mZ^2$,  let $f_{\ell_1 , \ell_2} \colon \mS^{15} \to B\oG_6$ be a map representing the homotopy class $\ell_1 \cdot [f] + \ell_2 \cdot [g] \in \pi_{15}(B\oG_6)$. 
 Let $M$ be the given manifold of dimension 15 with a trivial subbundle $\e_6 \subset TM$, and let  $\beta \colon M \to \mS^{15}$ be the map obtained by collapsing onto the top cell of $M$.  Then    the classifying maps $f_{\ell_1 , \ell_2} \circ \beta \colon M \to B\oG_6$ are pairwise non-homotopic, as they have distinct rigid secondary classes. Let $\F_{\ell_1 , \ell_2}$ be a codimension 6 foliation on $M$ constructed using the Thurston Realization Theorem, then 
 $\F_{\ell_1 , \ell_2}$ and $\F_{\ell_1' , \ell'_2}$ are not homotopic foliations for $(\ell_1 , \ell_2) \ne (\ell_12' , \ell_2')$. On the other hand, the foliations  $\F_{\ell_1 , \ell_2}$ on the fixed manifold $M$ have homotopic normal bundles, and so have homotopic tangent bundles.

 Again, for codimension $q=6$, the class $\Delta(y_2 y_4 c_2^3) \in H^{22}(B\oG_6 ; \mR)$ is spherically supported, so we can repeat the   procedure as above to obtain an infinite family $\{ \F_{\ell} \mid \ell \in \mZ\}$ of codimension 6 foliations on a manifold $M$ of dimension 22, which are pairwise non-homotopic, but have homotopic tangent bundles.
 
 For all even codimensions $q \geq 8$, we can repeat the above methods to obtain infinite families of foliations on compact manifolds $M$ of dimensions $n=2q+3, 2q+10, \ldots$ according to the degrees of the elements in $\ocR_q^s$, such that  the foliations are pairwise non-homotopic, but have homotopic tangent bundles.

\end{document}